\documentclass{article}

\usepackage{amssymb}
\usepackage{amsmath}
\usepackage[all]{xypic}
\usepackage[frenchb]{babel}
\usepackage{a4wide}
\usepackage{enumerate,fancyhdr}

\newcommand{\pp}{{\bf P}}
\newcommand{\qq}{{\bf Q}}

\newcommand{\oo}{\mathcal{O}}

\newcommand{\ra}{\rightarrow}

\newcommand{\mult}{{\rm mult}}

\newcommand{\Stab}{{\rm Stab}}
\newcommand{\codim}{{\rm codim}}

\newcommand{\C}{\mathbb{C}}
\newcommand{\card}{\textnormal{Card}}

\newcommand{\e}{\varepsilon}
\newcommand{\G}{\mathbb{G}}
\renewcommand{\L}{\mathcal{L}}
\newcommand{\demo}{\noindent \textit{D{\'e}monstration} : }

\newcommand{\cd}{\textnormal{cd}}
\newcommand{\stab}{\textnormal{Stab}}

\renewcommand{\phi}{\varphi}

\newcommand{\EV}{\mathbb{E}}
\newcommand{\N}{\mathbb{N}}

\newcommand{\tors}{\textnormal{tors}}

\newcounter{ndefinition}[section]
\newcommand{\defi}{\addtocounter{ndefinition}{1}{\noindent \textbf{D{\'e}finition \thendefinition.\ }}}
\newcounter{nrem}
\newcommand{\remarque}{\addtocounter{nrem}{1}{\noindent \textbf{Remarque \thenrem.\ }}}

\newtheorem{theo}{Th\'eor\`eme}
\newtheorem{lemm}[theo]{Lemme}

\newtheorem{prop}[theo]{Proposition}

\begin{document}

\title{Lemmes de multiplicit\'es et constante de Seshadri}

\author{Michael Nakamaye\footnote{Department of Mathematics and Statistics, University of New Mexico, Albuquerque, New Mexico, 87131, USA, nakamaye@math.unm.edu}\ et Nicolas Ratazzi \footnote{Universit\'e Paris-Sud XI, Batiment 425, 91405 Orsay Cedex, FRANCE, nicolas.ratazzi@math.u-psud.fr}}

\date{}

\maketitle

\hrulefill

\bigskip

\noindent \textbf{R\'esum\'e : } On d\'emontre dans cet article un raffinement des lemmes de multiplcit\'es de Philippon, essentiellement dans le cas particulier o\`u l'on d\'erive dans toutes les directions. L'am\'elioration est rendue possible grace \`a un point de vue plus g\'eom\'etrique et notamment par l'apparition nouvelle dans ce contexte de la notion de constante de Seshadri.

\bigskip

\hrulefill

\noindent \textbf{Abstract: } We establish an improvement of Philippon's
zero estimates primarily in the multiplicity setting. The improvement is
made possible by a more geometric approach and in particular the use of
Seshadri constants.

\bigskip

\hrulefill

\bigskip

\tableofcontents


\section{Introduction}

\subsection{Introduction}

\noindent On prouve dans cet article deux r\'esultats concernant les lemmes de multiplicit\'es. Avant d'aborder les \'enonc\'es proprement dit rappelons rapidement ce qu'est un lemme de multiplicit\'es et dans quel contexte on rencontre ce type d'\'enonc\'e : il s'agit de fournir des conditions suffisantes pour qu'une section d'un fibr\'e en droites $\L$ sur une compactification \'equivariante lisse $X$ d'un groupe alg\'ebrique commutatif $\G$ ne puisse pas s'annuler sur un sous-sch\'ema fini $\Sigma$ donn\'e sans \^etre identiquement nulle. Ce type d'information qui se traduit g\'en\'eralement par la pr\'esence d'un sous-groupe obstructeur dont on majore le degr\'e g\'eom\'etrique est d'usage particuli\`erement important en transcendance et en g\'eom\'etrie diophantienne et intervient \'egalement en g\'eom\'etrie alg\'ebrique complexe dans un contexte plus g\'en\'eral : partant d'une vari\'et\'e alg\'ebrique complexe projective et lisse $X$, munie d'un fibr\'e en droites ample $L$
sur $X$, et $\Sigma \subset X$ un sous--ensemble fini, un lemme de multiplicit\'es
pour les donn\'ees $X,L,\Sigma$ consiste \`a majorer la multiplicit\'e maximale que peut 
avoir une section non triviale $s \in H^0(X,L)$ le long de $\Sigma$.  Ce probl\`eme
est aussi difficile qu'important, m\^eme dans le cadre le plus simple.  
Par exemple, si $X = \pp^2$, $L = \oo_{\pp^2}(d)$, 
et $\Sigma$ est un ensemble de $m$
points g\'en\'eraux avec $m \geq 10$, on tombe  sur la conjecture de Nagata. 

\medskip

\noindent Deux aspects nous semblent importants dans un lemme de multiplicit\'es : le premier est d'obtenir une bonne condition num\'erique concernant le degr\'e du groupe ou de la vari\'et\'e obstructrice, le second est aussi de mieux comprendre g\'eom\'etriquement d'o\`u provient l'obstruction. De ce point de vue il nous semble int\'eressant de pouvoir identifier en terme g\'eom\'etrique les vari\'et\'es obstructrices. Dans cet article, nous nous int\'eressons aux lemmes de multiplicit\'es dans le cadre o\`u ils apparaissent en  transcendance. Dans ce contexte, le r\'esultat le plus important est d\^u \`a Philippon \cite{philippon} (cf. \'egalement \cite{nak1} pour une preuve plus g\'eom\'etrique, bas\'ee en partie sur celle de Philippon mais permettant \'egalement d'obtenir un r\'esultat plus pr\'ecis au niveau des constantes intervenant dans la majoration du degr\'e du sous-groupe obstructeur). En utilisant des outils d\'evelopp\'es en g\'eom\'etrie alg\'ebrique complexe nous donnons un raffinement, dans le cas particulier o\`u l'on d\'erive dans toutes les directions, des lemmes de z\'eros de Philippon \cite{philippon}. Par ailleurs, dans le cas de dimension $2$ nous obtenons \'egalement un r\'esulat plus pr\'ecis dans le cas o\`u l'on d\'erive le long d'une droite (\textit{cf.} le th\'eor\`eme \ref{t5} paragraphe \ref{observation}). Enfin utilisant une preuve plus g\'eom\'etrique  et notamment en utilisant la notion de constante de Seshadri nous donnons une information g\'eom\'etrique sur les sous-groupes obstructeurs, les reliant \`a la notion de sous-vari\'et\'e exceptionnelle de Seshadri.  En dimension deux le r\'esultat est particuli\`erement frappant~: (sous certaines conditions de taille sur l'ensemble de points $\Sigma$ consid\'er\'e fourni dans les donn\'ees) les sous-groupes obstructeurs sont pr\'ecis\'ement les vari\'et\'es exceptionnelle de Seshadri.

\medskip

\noindent Soient $\G/\C$ un groupe alg\'ebrique commutatif de dimension $d\geq 1$, $X$ une compactification de Serre de $\G$ (cf. \cite{serre}) (une telle compactification est notamment \'equivariante  et lisse). Si $U$ est une sous-vari\'et\'e de $\G$ nous noterons $\overline{U}$ son adh\'erence de Zariski dans $X$. Soient $\L$ un fibr\'e en droites ample sur $X$ et $\Gamma$ un sous-ensemble fini de $\G(\C)$ contenant $0$ et engendrant un groupe (n\'ecessairement de type fini) $\Gamma_{\mathbb{Z}}$. On pose pour tout entier $S\geq 1$
\[\Gamma(0)=\Gamma\ \text{ et }\ \Gamma(S)=\left\{\sum_{i=1}^Sx_i\ /\ \forall i\in\{1,\ldots,S\}, \ x_i\in\Gamma\right\}.\]

\noindent Par ailleurs, \'etant donn\'e un sous-espace vectoriel non nul $\EV$ de l'espace tangent \`a l'origine $T_0(\G)$ et \'etant donn\'e un sous-groupe alg\'ebrique $H$ de $\G$, nous noterons
\[c_H(\EV):=\codim(\EV\cap T_0(H), \EV)\text{ et } \alpha(S,\EV,\L)=\left(\frac{\deg_{\L}X}{|\Gamma(S)|}\right)^{\frac{1}{\dim \EV}}.\]
\noindent Dans le cas o\`u l'espace vectoriel consid\'er\'e $\EV$ est l'espace tangent $T_0(\G)$ entier, nous noterons simplement $\alpha(S,\L)$ plut\^ot que $\alpha(S,\EV,\L)$. Notons que si $D\geq 1$ est un entier, on a $\alpha(S,\EV,\L^{\otimes D})=D^\frac{\dim X}{\dim E}\alpha(S,\EV,\L)$.

\medskip

\noindent Nous pouvons maintenant \'enoncer le lemme de z\'eros de Philippon (dans la version pr\'ecis\'ee de Nakamaye)~:

\medskip

\noindent \textbf{Th\'eor\`eme}\label{t1}\footnote{Notons que l'\'enonc\'e classique de Philippon tel qu'on le trouve dans son article \cite{philippon} peut a priori sembler l\'eg\'erement moins g\'en\'eral puisqu'il suppose $S=a_1=...=a_{d-1}$. En fait il est imm\'ediat de voir  que l'\'enonc\'e ci-dessus se d\'eduit de celui obtenu avec $S=a_1=...=a_{d-1}$. }\textit{ Soient $S, a_1, \ldots, a_{d-1}$ des entiers strictement positifs ordonn\'es de fa\-\c{c}on d\'e\-crois\-sante et $T$ un entier strictement positif. Soient $D\geq 1$ un entier, $\EV$ un sous-espace vectoriel de l'espace tangent \`a l'origine $T_0(\G)$ et $\sigma$ une section non nulle de $H^0(X,\L^{\otimes D})$ s'annulant sur l'ensemble $\Gamma(S+a_1+\cdots+a_{d-1})$ le long de $\EV$ \`a un ordre au moins $dT+1$. Alors
il existe un sous-groupe alg\'ebrique $H$, diff\'erent de $\G$, tel que 
\[\max\{1,T\}^{c_H(\EV)}\card\left(\Gamma(a_{d-1})+H/H\right)\deg_{\L}(\overline{H})\leq D^{\dim X-\dim H}\deg_{\L}(X).\]
\noindent Le sous-groupe $H$ est appel\'e sous-groupe obstructeur.}

\medskip

\remarque Dans la pratique en transcendance, le fibr\'e en droites ample $\L$ est simplement la donn\'ee d'un plongement de la vari\'et\'e $X$ dans un espace projectif. C'est une donn\'ee du probl\`eme. Dans une preuve de transcendance on construit une section de ``degr\'e $D$" relativement \`a $\L$ (autrement dit un \'el\'ement de $H^0(X,\L^{\otimes D})$) nulle \`a un grand ordre $T$ sur un certain ensemble $\Sigma$ et un lemme de multiplicit\'es consiste essentiellement \`a obtenir une majoration non triviale de $T$  en fonction de $D$.

\medskip

\noindent Dans la suite nous ferons l'hypoth\`ese restrictive suivante qui justifie dans le titre le choix du terme \textit{lemme de multiplicit\'es} plut\^ot que  \textit{lem\-me de z\'eros}~: 

\medskip

\begin{tabular}{lr}
\textbf{Hypoth\`ese : } 	& $T\geq 1$.
\end{tabular}

\medskip

\noindent Sous cette hypoth\`ese, la conclusion du th\'eor\`eme pr\'ec\'edent peut \^etre reformul\'ee sous forme d'une dichotomie : 
\begin{enumerate}
\item Si $T\leq \alpha(S,\EV,\L^{\otimes D})$ alors le groupe trivial $\{0\}$ est obstructeur~:
\[T^{c_{\{0\}}(\EV)}|\Gamma(S)|\leq D^{\dim X}\deg_{\L}(X).\]
\item Sinon, il existe un sous-groupe alg\'ebrique $H$, diff\'erent de $\G$, tel que 
\[T^{c_H(\EV)}\card\left(\Gamma(a_{d-1})+H/H\right)\deg_{\L}(\overline{H})\leq D^{\dim X-\dim H}\deg_{\L}(X).\]
\end{enumerate}
\noindent Notons que dans cette seconde partie de l'alternative rien ne dit que le groupe $H$ ne peut pas \^etre le groupe trivial. Une des choses que nous faisons \'egalement dans cette article est de montrer que sous certaines conditions, on peut assurer que $H$ est non trivial.

\medskip

\noindent Cette formulation \`a l'avantage de bien indiquer que dans tous les lemmes de multiplicit\'es ($T\geq 1$), on peut en fait supposer que 
\[T>\alpha(S,\EV,\L^{\otimes D}),\] 
\noindent ce que l'on fera d\'esormais ; l'autre cas \'etant en fait trivial (rappelons que la philosophie des lemmes de multiplicit\'es est la suivante : partant d'une section qui s'annule \`a un grand ordre en un certain nombre de points suffisamment bien r\'epartis, on montre que l'ordre ne peut en fait pas \^etre trop grand. Si l'on part d'un ordre d'annulation $T\geq 1$ qui est d\'eja petit, il n'y a, du point de vue des lemmes de multiplicit\'es, rien de plus \`a dire). 

\medskip

\noindent L'objectif de cet article est de donner des informations suppl\'ementaires sur le sous-groupe obstructeur dans le cas particulier o\`u l'on fixe l'espace des d\'erivations \`a l'espace tangent $T_0(\G)$ tout entier, plut\^ot qu'\`a un sous-espace vectoriel quelconque de celui-ci. Dans ce cas nous montrons que l'on peut l\'eg\`erement affaiblir l'hypoth\`ese sur l'ordre d'annulation de la section $\sigma$ consid\'er\'ee tout en obtenant une information plus pr\'ecise sur le groupe obstructeur. Par ailleurs, dans le cas o\`u le groupe $\Gamma_{\mathbb{Z}}$ engendr\'e par $\Gamma$ n'est pas un groupe de torsion, on donne une condition num\'erique sur les param\`etres $a_i$ assurant que le groupe  obstructeur $H$ est diff\'erent du groupe trivial r\'eduit \`a l'\'el\'ement neutre. En particularisant la situation en dimension $2$ nous montrons qu'une condition num\'erique plus faible est suffisante dans ce cas. Enfin, toujours dans le cas particulier de la dimension $2$, nous expliquons au paragraphe \ref{observation} (th\'eor\`eme \ref{t5}) comment faire \'egalement fonctionner notre approche dans le cas d'un espace de d\'erivations de dimension $1$ (inclus dans l'espace tangent de dimension $2$). 

\subsection{\'Enonc\'e des r\'esultats}

\defi Soient $X/\C$ une vari\'et\'e alg\'ebrique, irr\'eductible et projective, et $\L$ un fibr\'e en droites ample sur $X$. Pour tout ensemble fini $\Gamma\subset X$, la \textit{constante de Seshadri} $\e(\Gamma,\L)$ est d\'efinie par 
\[\e(\Gamma,\L)=\inf_{C\subset X,\ C\cap\Gamma\not=\emptyset}\frac{\L\cdot C}{\sum_{x\in\Gamma}\mult_x C},\]
\noindent o\`u $C$ d\'ecrit les courbes irr\'eductibles de $X$ contenant au moins un point de $\Gamma$.

\medskip

\noindent On note que $\e(\Gamma,\L^{\otimes D})=D\e(\Gamma,\L)$. Par ailleurs rappelons (cf. par exemple \cite{lazarsfeld} p. 271 Proposition 5.1.9) que pour toute sous-vari\'et\'e irr\'eductible $V$ de $X$, de dimension non nulle et contenant au moins un point de $\Gamma$, on a 
\begin{equation}\label{varex}
\e(\Gamma,\L)\leq\left(\frac{\L^{\dim V}\cdot V}{\sum_{x\in\Gamma}\mult_x V}\right)^{\frac{1}{\dim V}}.
\end{equation}
\noindent Notamment en appliquant ceci avec $V=X$ on obtient : 
\begin{equation*}
\e(\Gamma,\L)\leq\left(\frac{\deg_{\L}X}{\sum_{x\in\Gamma}\mult_x X}\right)^{\frac{1}{\dim X}}.
\end{equation*}
\noindent Il n'est pas du tout clair a priori qu'il existe une vari\'et\'e telle que l'in\'egalit\'e dans (\ref{varex}) soit en fait une \'egalit\'e. N\'eanmoins le th\'eor\`eme de Campana-Peternell \cite{cp} sur le crit\`ere de Nakai pour les $\mathbb{R}$-diviseurs implique qu'il existe toujours une vari\'et\'e $V$ r\'ealisant l'\'egalit\'e (cf. \cite{lazarsfeld} p. 271 Proposition 5.1.9). Nous renvoyons au livre \cite{lazarsfeld} pour des faits g\'en\'eraux sur la constante de Seshadri. Ce que nous appelons dans cet article constante de Seshadri, est \'egalement appel\'e \textit{multiple point Seshadri constant} dans la litt\'erature \textit{cf.} par exemple Bauer \cite{bauer}.

\medskip

\defi Une vari\'et\'e $V$ r\'ealisant l'\'egalit\'e dans (\ref{varex}) est appel\'ee \textit{vari\'et\'e exceptionnelle de Seshadri de $\L$ relativement \`a $\Gamma$}. 

\medskip

\remarque \label{re1}Notons que si $V$ est une vari\'et\'e exceptionnelle de Seshadri pour $\L$, alors, par homog\'en\'eit\'e de $\e$, c'est \'egalement une vari\'et\'e exceptionnelle de Seshadri pour $\L^{\otimes D}$.

\medskip

\noindent \textbf{Notations : } \'etant donn\'ee une 
sous-vari\'et\'e $V$ de $X$, nous noterons 
\[\stab(V)=\left\{x\in\G\ / \ x+V=V\right\}\] 
\noindent le \textit{stabilisateur de $V$}. C'est un sous-groupe alg\'ebrique de $\G$ de composante neutre (\textit{i.e.} la composante connexe contenant l'identit\'e) $\stab(V)^0$. Enfin, \'etant donn\'ee une sous-vari\'et\'e $U$ de $\G$ nous noterons $\overline{U}$ son adh\'erence de Zariski dans la compactification $X$.

\noindent Par ailleurs nous noterons $u$ la solution r\'eelle strictement comprise entre $0$ et $1$ de l'\'equation $x^{d-1}(1+x)=1$ avec $d=\dim X$.

\medskip

\noindent Rappelons nos notations :  $\G/\C$ est un groupe alg\'ebrique commutatif de dimension $d\geq 1$, $X$ une compactification de Serre de $\G$  et $\Gamma$ un sous-ensemble fini de $\G(\C)$ contenant $0$ et engendrant un groupe  $\Gamma_{\mathbb{Z}}$.

\medskip

\begin{theo}\label{t2}Soient $T$ et $D$ deux entiers strictements positifs, $\L$ un fibr\'e en droites ample, $S=a_0, a_1, \ldots, a_{d-1}$ des entiers strictement positifs ordonn\'es de fa\c{c}on d\'ecroissante et $V$ une vari\'et\'e exceptionnelle de Seshadri de $\L$ relativement \`a $\Gamma(S)$. On suppose que $T>\alpha(S,\L^{\otimes D})$. Soit alors $\sigma$ une section non nulle  de $H^0(X,\L^{\otimes D})$ s'annulant sur l'ensemble $\Gamma(S+a_1+\cdots+a_{\codim V})$ le long de $T_0(\G)$ \`a un ordre sup\'erieur ou \'egal \`a $(u+\codim V)T$ . 
\begin{enumerate}
\item Il existe une vari\'et\'e $W$ de dimension comprise entre $\dim V$ et $\dim X-1$ contenant un translat\'e de $V$ par un point de $\Gamma(a_1+\cdots+a_{\codim V-1})$ telle qu'en posant $H=\stab(W\cap \G)^0$ on a 

\[T^{\codim H}\card\left(\Gamma(a_{\codim V})+H/H\right)\deg_{\L}(\overline{H})\leq\left(\deg_{\L}(X)\right)D^{\codim H}.\]
\item Si on suppose de plus que $\Gamma_{\mathbb{Z}}$ n'est pas de torsion, alors en choisissant 
\[\forall i\in\{1,\ldots,\codim V\},\ \ |\Gamma(a_i)|>|\Gamma_{\mathbb{Z},\textnormal{tors}}|(\deg_{\L}X) |\Gamma(S)|^{\frac{\codim V}{d}},\]
\noindent on a de plus~: le sous-groupe strict $H$ obtenu au point 1. est non nul. 

\end{enumerate}
\end{theo}

\medskip

\noindent Ce th\'eor\`eme \ref{t2} appelle un certain nombre de commentaires et remarques : 

\medskip

\begin{enumerate}
\item Pour le point 1. du th\'eor\`eme, plus le nombre 
\[\sup\{\dim V\ / \ \text{$V$ vari\'et\'e exceptionnelle de $\L$ relativement \`a $\Gamma(S)$}\}\]
\noindent est grand et meilleur est le r\'esultat . 
\item Au vu de la preuve (cf. le lemme \ref{rajout}), on peut en fait remplacer $u$ par une constante un peu meilleure : l'unique solution dans $]0,1[$ de l'\'equation $x^{d-1}(x+\codim V)=1$ o\`u $V$ est la vari\'et\'e exceptionnelle de Seshadri intervenant dans l'\'enonc\'e du th\'eor\`eme.
\item Dans le point 2. du th\'eor\`eme il est important de noter que dans les hypoth\`eses faites sur les nombres $a_i$, on ne fait intervenir que le degr\'e de $X$ relativement \`a $\L$ et non pas relativement \`a $\L^{\otimes D}$ : le choix des $a_i$ ne d\'epend pas de $D$. 
\end{enumerate}

\medskip

\noindent Les apports de cet \'enonc\'e nous semblent \^etre les suivants~:

\begin{enumerate}
\item Si l'on se restreint au cas des surfaces, il permet de montrer que toute courbe exceptionnelle de Seshadri donne directement naissance, par passage au stabilisateur, \`a un sous-groupe obstructeur ; et est elle m\^eme un sous-groupe obstructeur si les hypoth\`eses du point 2. sont satisfaites. En dimension quelconque, la vari\'et\'e que l'on construit naturellement, avant de passer au stabilisateur pour obtenir un groupe obstructeur, contient un translat\'e de vari\'et\'e exceptionnelle de Seshadri.
\item On demande une condition d'annulation plus faible que le classique $dT+1$ : notre condition d'annulation est au pire de la forme $(d-1+u)T$ avec $0< u < 1$. Ceci est d\^u \`a l'introduction de la notion de constante et de vari\'et\'e de Seshadri et \'evidemment n'est rendu possible que par l'hypoth\`ese faite sur l'espace des d\'erivations : \^etre tout l'espace tangent $T_0(\G)$.
\item Dans le cas o\`u l'on part d'un ensemble engendrant un groupe non de torsion, on prouve une version plus forte des lemmes de multiplicit\'es en partant d'une section nulle sur un ensemble plus petit que $\Gamma(dS)$, avec un ordre moins grand que classiquement, et en assurant n\'eanmoins l'existence d'un sous-groupe obstructeur strict \textit{non nul} dans ce cas. 
\item Si $\Gamma$ n'est pas contenu dans $\G_{\rm tors}$, il suffit pour obtenir les in\'egalit\'es dans le point 2. du Th\'eor\`eme \ref{t2} de choisir les $a_i$ tel que
\[
|\Gamma(a_i)| = o(S).
\]
Autrement dit, au lieu de supposer que $\sigma$ s'annule sur $\Gamma(dS)$, il suffit dans ce cas de ne supposer l'annulation que sur
$\Gamma(S + o(S))$.  En fonction de $S$, ceci est le meilleur r\'esultat
possible.  
\item La raison profonde pour laquelle le Th\'eor\`eme \ref{t2} ne se
g\'en\'eralise pas bien au cas d'un sous-espace quelconque $\EV \subset 
T_0(\G)$, est que la notion d'ordre d'annulation le long de $\EV$
s'adapte mal \`a une interpr\'etation g\'eom\'etrique. Autrement dit, alors
que la constante de Seshadri $\epsilon(x,\L)$ est une fonction homog\`ene
de $\L$, l'analogue qui mesure la positivit\'e de $\L$ le long de $\EV$ ne l'est plus.  
\end{enumerate}

\medskip

\noindent Ceci \'etant dit, il convient de temp\'erer notre r\'esultat : il est en g\'en\'eral tr\`es difficile de calculer la constante de Seshadri et \`a plus forte raison de donner des informations sur les vari\'et\'es exceptionnelles de Seshadri (dont on sait tout de m\^eme qu'elles sont au moins de dimension $1$).

\medskip

\noindent Dans le cas des surfaces nous pouvons, pour le point 2. du th\'eor\`eme \ref{t2} pr\'ec\'edent, l\'eg\`erement affaiblir les hypoth\`eses~:

\begin{theo}\label{t4} On suppose que $\dim X=2$ et que $\Gamma_{\mathbb{Z}}$ n'est pas de torsion. Soient $D\geq 1$ un entier, et $S$, $a$ deux entiers strictement positifs tels que $|\Gamma(S)|\geq |\Gamma(a)|\geq|\Gamma_{\mathbb{Z},\textnormal{tors}}|\cdot|\Gamma(S)|^{\frac{1}{2}}$. Soient $T>\alpha(S,\L^{\otimes D})$ un entier et $\sigma \in H^0(X,\L^{\otimes D})$ une section non nulle s'annulant \`a un ordre au moins $\frac{1}{2}\left(1+\sqrt{5}\right)T$ sur $\Gamma(S+a)$. Toute sous-vari\'et\'e exceptionnelle de Seshadri de $\L$ relativement \`a $\Gamma(S)$ est l'adh\'erence d'un translat\'e de sous-groupe alg\'ebrique de dimension $1$. De plus, si $E \subset X$  est l'une de ces courbes, alors $\sigma$ s'annule le long de $\Gamma(a) + E$ \`a un ordre au moins $T$
 et l'on a
\[ T\cdot\card(\Gamma(a) + E/E)\deg_{\L}(E) \leq D\deg_{\L}(X).\]
\end{theo}

\medskip

\remarque Ce dernier r\'esultat n'est pas une cons\'equence du th\'eor\`eme \ref{t2} pr\'ec\'edent dans le cas des surfaces car l'hypoth\`ese faite sur la taille de $\Gamma(a)$ est plus faible que celle du th\'eor\`eme \ref{t2}~: elle ne fait pas intervenir le degr\'e $\deg_{\L}(X)$. La raison profonde sous-jacente est qu'en dimension $2$ la vari\'et\'e obstructrice est n\'ecessairement une hypersurface.

\medskip

\remarque On explique au paragraphe \ref{observation} (\textit{cf.} notamment le th\'eor\`eme \ref{t5}) comment traiter en dimension $2$ le cas d'un sous-espace $\EV$ de dimension $1$.

\medskip

\noindent Nous commen\c{c}ons par rappeler la notion de d\'erivation de sections que nous utiliserons, qui est une notion de d\'erivation locale. Nous donnons ensuite la preuve du premier th\'eor\`eme en insistant particuli\`erement sur la proposition \ref{l1} qui est le point cl\'e r\'eellement nouveau. Nous r\'eutilisons ensuite cette proposition dans la preuve du second th\'eor\`eme.

\medskip

\noindent \textbf{{\em Remerciements}} C'est un plaisir pour les deux auteurs de remercier les universit\'e
de Paris VII et de Paris-Sud Orsay qui ont cordialement
re\c{c}u le premier auteur pendant deux mois, nous permettant ainsi d'accomplir
ce travail ensemble.

\section{Les d\'erivations et le th\'eor\`eme de B\'ezout}

\noindent Nous utiliserons dans la suite une notion de d\'erivation locale. C'est cette notion que l'on trouve par exemple introduite dans \cite{szpiro} qui est utilis\'ee dans l'article \cite{nak1} et qui permet de gagner au niveau des constantes intervenant dans l'in\'egalit\'e de majoration du degr\'e du groupe obstructeur $H$. Nous \'enon\c{c}ons ensuite le th\'eor\`eme de B\'ezout que nous utiliserons dans la preuve du th\'eor\`eme \ref{t2}.

\subsection{D\'erivations}

\noindent Soient $X$ une vari\'et\'e alg\'ebrique munie d'un fibr\'e en droites $\L$, $\sigma$ une section non nulle de $H^0(X,\L)$ et $D$ un op\'erateur diff\'erentiel  d'ordre $1$ sur $X$, \textit{i.e.} tel que 
\[\forall f,g\in\oo_X,\ \ \ D(fg)=D(f)g+fD(g).\]
\noindent La d\'eriv\'ee $D(\sigma)$ de $\sigma$ selon $D$ n'est pas toujours une section bien d\'efinie de $H^0(X,\L)$. Par contre, pour toute sous-vari\'et\'e $Y$ contenue dans le lieu des z\'eros $Z(\sigma)$ de $\sigma$, il est possible de d\'efinir $D(\sigma)$ comme un \'el\'ement de $H^0(Y,\L_{|Y})$. Consid\'erons pour cela une trivialisation $\{U_i\}_{i\in I}$ de $\L$, avec des fonctions de transitions $\phi_{ij}$ sur $U_i\cap U_j$, telle que la section $\sigma$ est donn\'ee sur $U_i$ par des fonctions r\'eguli\`eres $f_i$. On a donc
\[\forall i,j\in I,\ \ \  f_j=\phi_{ij}f_i.\]
\noindent En composant par $D$ on en d\'eduit que sur $U_i\cap U_j$ on a
\[\forall i,j\in I,\ \ \  D(f_j)=D(\phi_{ij})f_i+\phi_{ij}D(f_i).\]
\noindent On voit ainsi que la condition de recollement n'est pas forc\'ement assur\'ee globalement. Par contre pour toute sous-vari\'et\'e $Y\subset Z(\sigma)$, on obtient bien $D(f_j)|Y=\phi_{ij}D(f_i)|Y$ et ceci nous indique par recollement que l'on a une section bien d\'efinie $D(\sigma)\in H^0(Y,\L_{|Y})$.

\medskip

\noindent Ceci permet par it\'eration de d\'efinir une section $D(\sigma)$ pour tout op\'erateur diff\'erentiel $D$. La section ainsi construite est bien d\'efinie sur toute sous-vari\'et\'e $V$ de $X$ telle que $D_1(\sigma)_{|V}\equiv 0$ pour tout op\'erateur $D_1\ll D$ avec des notations \'evidentes (on prend une base de l'espace des op\'erateurs et on dit que $D_1\ll D$ si les seuls \'el\'ements de la base apparaissant dans $D_1$ apparaissent aussi dans $D$ et si l'ordre auquel ils apparaissent dans $D_1$ est inf\'erieur ou \'egal \`a celui auquel ils apparaissent dans $D$ et strictement inf\'erieur pour au moins l'un d'entre eux). 

\subsection{Th\'eor\`eme de B\'ezout}

\noindent Nous rappelons ici le th\'eor\`eme de B\'ezout que nous utiliserons, dans une version que l'on trouve dans Fulton \cite{fulton} corollaire du paragraphe 12.4.8. Soient $D_1,\ldots,D_d$ des diviseurs effectifs sur une vari\'et\'e $V$ de dimension $d$. Si $x$ est un point isol\'e de l'intersection de $D_1,\ldots,D_d$ sur $V$, alors on note 
\[i(x,D_1\cdot\ldots\cdot D_d; V)\]
\noindent la multiplicit\'e d'intersection de $x$, \textit{i.e.} le coefficient de $[x]$ dans la classe $D_1\cdot\ldots\cdot D_d\cdot V$.

\medskip

\begin{prop}\label{bezout}Soient $D_1,\ldots,D_d$ des diviseurs effectifs sur une vari\'et\'e alg\'ebrique $V$ de dimension $d$. Si $x$ est un point isol\'e de l'intersection de $D_1,\ldots,D_d$ sur $V$, alors
\[\mult_x(V)\prod_{k=1}^d\mult_x(D_k)\leq i(x,D_1\cdot\ldots\cdot D_d ; V).\]
\end{prop}

\section{Preuve du point 1. du th\'eor\`eme \ref{t2}}

\noindent L'essentiel de l'argument est bas\'e sur la m\'ethode maintenant classique de Philippon \cite{philippon} avec cependant un ingr\'edient nouveau qui est l'introduction de la notion de constante de Seshadri. Le point v\'eritablement nouveau \'etant contenu dans la proposition \ref{l1}. Sch\'ematiquement la preuve est la suivante~: on fixe une vari\'et\'e exceptionnelle de Seshadri $V$ de dimension $\dim V\geq 1$ et de codimension $\cd(V)$. En posant $V=V_{\cd(V)}$ et en utilisant que $0\in \Gamma$, on construit une suite emboit\'ee de sous-ensembles alg\'ebriques  de $X$, 
\[V_{\cd(V)}\subset V_{\cd{V}-1}\subset\ldots\subset V_{0}=Z(\sigma)\subsetneq X\]
\noindent telle que chaque inclusion est avec multiplicit\'e $T$ et telle que 
\[\forall k\in\left\{0,\ldots,\cd(V)-1\right\}\ \ \ \ \Gamma(a_{\cd(V)-k})+V_{k+1}\subset V_k\]
\noindent avec l\`a encore des inclusions avec multiplicit\'e $T$.
\noindent  Par le principe des tiroirs et en oubliant les multiplicit\'es, on en d\'eduit que au moins deux des $V_i$ sont de m\^eme dimension, \textit{i.e.} il existe $\cd(V)-1\geq k\geq 0$ tel que $\dim V_{k+1}=\dim V_{k}$. On note $i$ le plus petit entier v\'erifiant ceci. Le reste de la preuve est ensuite exactement calqu\'e sur la m\'ethode de Philippon, revisit\'ee par Nakamaye : on note $W_i$ une composante irr\'eductible de $V_{i+1}$ de dimension maximale. C'est aussi une composante irr\'eductible isol\'ee de $V_i$ et mieux, $\Gamma(a_{\cd(V)-i})+W_i$ est une somme de composantes irr\'eductibles isol\'ees de $V_i$. Utilisant que l'inclusion est avec multiplicit\'e $T$ on construit, en coupant par des diviseurs num\'eriquement \'equivalents \`a $\L^{\otimes D}$ un cycle dont le support contient  $\Gamma(a_{\cd(V)-i})+W_i$ avec multiplicit\'e $T$. Le th\'eor\`eme de B\'ezout nous donne alors l'in\'egalit\'e
\[T^{\codim W_i}\deg_{\L^{\otimes D}}\left(\Gamma(a_{\cd(V)-i})+W_i\right)\leq\deg_{L^{\otimes D}} X.\]
\noindent Il ne reste plus ensuite qu'\`a obtenir le groupe obstructeur. Pour cela on consid\`ere le stabilisateur de $W_i$ :  c'est une composante isol\'ee de translat\'es de $V_{i}$ avec multiplicit\'e $T$, donc la m\^eme application du th\'eor\`eme de B\'ezout permet de conclure. Le point nouveau ici est que l'on est assur\'e dans le d\'ebut de la preuve de d\'emarrer la filtration avec un ensemble $V_{\cd(V)}$ de dimension $\dim V$, alors que dans les preuves classiques la seule information connue \'etait que l'ensemble $V_{\dim X}$ est non-vide. C'est cette nouveaut\'e qui permet de gagner sur les hypoth\`eses concernant l'ordre d'annulation de la section consid\'er\'ee $\sigma$. Le reste de l'argument est inchang\'e par rapport aux preuves de Philippon et Nakamaye.

\begin{lemm}\label{l0} Soient $X/\mathbb{C}$ une vari\'et\'e alg\'ebrique projective munie d'un fibr\'e en droites ample $\L$ et $\Gamma$ un ensemble fini de points de $X$. L'in\'egalit\'e
\[\e(\Gamma,\L)\leq \left(\frac{\L^d\cdot V}{\sum_{x\in\Gamma}\mult_x V}\right)^{\frac{1}{d}}\]
\noindent est valable pour toute sous-vari\'et\'e $V$ de pure dimension $d\geq 1$ de $X$, irr\'eductible ou non.
\end{lemm}
\demo On sait d\'ej\`a que le r\'esultat est vrai si $V$ est irr\'eductible. On en d\'eduit imm\'edia\-tement le cas g\'en\'eral en utilisant l'in\'egalit\'e
\[\inf\left\{\frac{a}{c},\frac{b}{d}\right\}\leq \frac{a+b}{c+d}\]
\noindent valable pour tout entier non nul $a,b,c,d$.\hfill$\Box$

\medskip

\noindent Soient $D,T,S,a_1,\cdots,a_{\cd(V)}$ des entiers strictement positifs, $\L$ un fibr\'e en droites ample, $V$ une vari\'et\'e exceptionnelle de Seshadri de $\L$ relativement \`a $\Gamma(S)$ et $\sigma$ une section non nulle de $H^0(X,L^{\otimes D})$. On pose 
\[V_0=Z(\sigma),\ \ \ V_{\cd(V)}=V,\ \ \ \Gamma_{\cd(V)}=\Gamma\left(a_{1}+\cdots+a_{\cd(V)}\right)\]
\noindent et pour tout entier $k$ tel que $\cd(V)-1\geq k\geq 1$, on pose
\[\Gamma_k=\Gamma\left(a_{\cd(V)+1-k}+\cdots+a_{\cd(V)}\right)\text{ et } V_k=\left\{x\in X\ / \ \mult_{x+\Gamma_k}\sigma\geq kT+1\right\}.\]

\begin{prop}\label{l1} On suppose que la section $\sigma$ s'annule sur $\Gamma(S)+\Gamma_{\cd(V)}$ \`a un ordre au moins $(\cd(V)+u)T$ o\`u $u\in ]0,1[$ est un r\'eel tel qu'on a l'in\'egalit\'e $uT>\e(\Gamma(S),\L^{\otimes D})$. Dans ce cas on a
\[\Gamma(a_1)+V\subset V_{\cd(V)-1}.\]
\noindent De plus il s'agit d'une inclusion avec multiplicit\'e au moins $T$, \textit{i.e.}, 
\[\Gamma(a_1)+V\subset\hspace{-.3cm} \bigcap_{\overset{D\textnormal{op. diff.}}{\textnormal{d'ordre au plus T}}}\left\{x\in X\ / \ \mult_{x+\Gamma_{\cd(V)-1}}D(\sigma)\geq (\cd(V)-1)T+1\right\}.\]
\end{prop}
\demo Soit $s\in H^0(X,\L^{\otimes D})$ une section nulle \`a un ordre au moins $uT$ sur $\Gamma(S)+\Gamma_{\cd(V)}$. Pour tout point $g\in \Gamma_{\cd(V)}$, la section $s$ est nulle \`a un ordre au moins $uT$ sur $\Gamma(S)+g$, donc la section $t^*_g(s)\in H^0(X,t^*_g\L^{\otimes D})$ est nulle \`a un ordre au moins $uT$ sur $\Gamma(S)$. Supposons par l'absurde que cette derni\`ere section n'est pas identiquement nulle sur la vari\'et\'e exceptionnelle $V$.

\medskip

\noindent Consid\'erons tout d'abord le cas plus simple o\`u $\dim V=1$. Dans ce cas il suffit simplement d'appliquer le th\'eor\`eme de B\'ezout \`a l'intersection 
$Z(t_g^\ast(s)) \cap V$ pour aboutir \`a l'in\'egalit\'e 
\begin{eqnarray}
D\deg_{\L}(V)	& =	&	\deg_{\L^{\otimes D}}(V)\geq Z\left(t_g^\ast(s)\right)\cap V\\
				& \geq 	& \sum_{x\in \Gamma(S)}uT\mult_x(V)=uT\frac{\deg_{\L}(V)}{\e(\Gamma(S),\L)}.
\end{eqnarray}
\noindent En simplifiant on obtient donc
\begin{equation}\label{f0}
uT\leq D\e(\Gamma(S),\L)=\e\left(\Gamma(S),\L^{\otimes D}\right).
\end{equation}
\noindent Ceci contredit l'hypoth\`ese faite sur $u$ et assure donc que la section $s$ est nulle sur $g+V$ et donc sur $\Gamma_{\cd(V)}+V$.

\medskip

\noindent Dans le cas g\'en\'eral d'une vari\'et\'e $V$ de dimension $\dim V\geq 2$ on veut aboutir \`a la m\^eme in\'egalit\'e pour conclure de la m\^eme fa\c{c}on. Malheureusement dans ce cas, une simple application du th\'eor\`eme de B\'ezout ne suffit plus et les choses se compliquent un peu. Pr\'ecis\'ement nous allons utiliser une seconde d\'efinition de la constante de Seshadri (voir par exemple [Laz] Definition 5.1.1), \'equivalente \`a celle que nous avons donn\'ee~: soit $\pi: \tilde{X}\ra X$ l'\'eclatement de $X$ le long de $\Gamma(S)$, de diviseur exceptionnel
$E$.  La constante de Seshadri de $\L$ relativement \`a $\Gamma(S)$ est donn\'ee par la formule 
\[\e(\Gamma(S),\L)=\sup\left\{\lambda\in\mathbb{R}\ / \ \pi^*\L-\lambda E\text{ nef}\ \right\}.\]

\medskip

\noindent \textbf{Fait : }pour tout r\'eel $0<\alpha < \e(\Gamma(S),\L)$ le 
fibr\'e $\pi^\ast(\L)(-\alpha E)$ est ample.

Notons d'abord que
$\e(\Gamma(S),\L) > 0$.  En effet si $N>0$ est un entier
tel que $\L^{\otimes N}$ est tr\`es ample alors $\e(x,\L^{\otimes N}) \geq 1$
pour tout $x \in X$ et par cons\'equent $\pi^\ast \L^{\otimes N}-\e E$
est ample, donc nef, pour tout $0 < \e < \frac{1}{|\Gamma(S)|}\leq 1$.

Soit maintenant $0 < \alpha < \e(\Gamma(S),\L)$.
Par d\'efinition de la constante de Seshadri $\e(\Gamma(S),\L)$ de $\L$ relativement \`a $\Gamma(S)$, le fibr\'e $\pi^\ast \L - \alpha
E$ est nef.  Supposons par l'absurde qu'il n'est pas ample. Le th\'eor\`eme de
Campana-Peternell \cite{cp} sur le crit\`ere de Nakai pour les $\mathbb{R}$-diviseurs nous assure dans ce cas qu'il existe une sous-vari\'et\'e $Y$ de $X$ de dimension
 $\dim Y\geq 1$ telle que
\[\deg_{\pi^\ast(\L)(-\alpha E)}(\tilde{Y}) = 0,\]
o\`u $\tilde{Y}$ d\'esigne la transform\'ee stricte de $Y$ dans l'\'eclatement $\tilde{X}$. 
Puisque
\[\deg_{\pi^\ast(\L)(-\alpha E)}(\tilde{Y}) = \deg_\L(Y) - \sum_{x \in \Gamma(S)}
\alpha^{\dim(Y)}\mult_x(Y)\]
on en d\'eduit que si $\lambda$ est un r\'eel v\'erifiant $\alpha < \lambda < \e(\Gamma(S),\L)$, alors le fibr\'e $\pi^\ast(\L)(-\lambda E)$ n'est pas nef. Ceci contredit la
d\'efinition de la constante de Seshadri $\e(\Gamma(S),\L)$ et prouve donc le fait.

\medskip

\noindent Supposons maintenant que $\alpha < \e(\Gamma(S),\L)$ est un rationnel
strictement positif. On sait par le fait pr\'ec\'edent que
le fibr\'e $\pi^\ast(\L)(-\alpha E)$ est ample.
Soient alors $N>0$ un entier (d\'ependant de $\alpha$) tel que le fibr\'e $\pi^\ast \L^{\otimes N} - N \alpha E$
est tr\`es ample et $E_1,\ldots, E_{\dim(V)-1}$ des diviseurs g\'en\'eraux du syst\`eme lin\'eaire
$|\pi^\ast \L^{\otimes N} - N \alpha E|$.  Ces diviseurs \'etant g\'en\'eraux
ils sont lisses et s'intersectent proprement.  Par ailleurs
le fibr\'e $\pi^\ast \L^{\otimes N} - N \alpha E$ \'etant tr\`es ample
et les diviseurs $E_1,\ldots, E_{\dim(V)-1}$ \'etant g\'en\'eraux, on a : pour toute
composante irr\'eductible $W \subset \pi^{-1}(Z(t_g^\ast(s))) \cap \tilde{V}$ 
l'intersection
\[E_1 \cap \ldots \cap E_{\dim(V)-1} \cap  W\]
est propre.
En particulier, l'intersection
\[\pi(E_1) \cap \ldots \cap \pi(E_{\dim(V)-1}) \cap Z(t_g^\ast(s)) \cap V\] 
est propre.  En posant $D_i = \frac{\pi(E_i)}{N}$ pour tout entier $i$ compris entre $1$ et $\dim V-1$ on a par construction
\[\forall x\in\Gamma(S),\ \ \mult_x(D_i)\geq \alpha\]
\noindent et de plus l'intersection
\[D_1 \cap \ldots \cap D_{\dim(V)-1} \cap Z(t_g^\ast(s)) \cap  V\]
est propre.

\medskip

\noindent En utilisant la proposition \ref{bezout} (th\'eor\`eme de B\'ezout) dans la vari\'et\'e $V$ et le
fait que $t_g^\ast(\L^{\otimes D})$ est num\'eriquement \'equivalent 
au fibr\'e en droites $\L^{\otimes D}$ on obtient, en faisant tendre $\alpha$ vers 
$\e(\Gamma(S),\L)$,
\[\sum_{x \in \Gamma(S)} \left(D\e(\Gamma(S),\L)\right)^{\dim(V)-1}\mult_x(V)uT\leq \deg_{\L^{\otimes D}}(V).\]
De ceci on d\'eduit l'in\'egalit\'e
\begin{equation}\label{f1}
uT\e(\Gamma(S),\L)^{\dim(V)-1} \leq \frac{\deg_\L(V)}{\sum_{x \in \Gamma(S)} \mult_x(V)}D.
\end{equation}
La vari\'et\'e $V$ \'etant une sous-vari\'et\'e exceptionnelle de Seshadri de $\L$ relativement \`a $\Gamma(S)$, nous avons
\begin{equation}\label{f2}
\e(\Gamma(S),\L) = \left(\frac{\deg_\L(V)}{\sum_{x \in \Gamma(S)} \mult_x(V)}\right)^{\frac{1}{\dim(V)}}. 
\end{equation}
En combinant les deux formules (\ref{f1}) et (\ref{f2}) pr\'ec\'edentes on conclut comme lorsque $V$ \'etait suppos\'ee de dimension $1$ que
\[uT \leq \e(\Gamma(S),\L)D=\e\left(\Gamma(S),\L^{\otimes D}\right).\]
\noindent Ceci contredit l'hypoth\`ese de la proposition et prouve donc par l'absurde que $t_g^\ast(s)|V = 0$. On en d\'eduit que la section $s$ est nulle sur $g+V$ et donc que $s$ est nulle sur $\Gamma_{\cd(V)}+V$.

\medskip

\noindent Consid\'erons maintenant la section  $\sigma\in H^0(X,\L^{\otimes D})$ nulle sur $\Gamma(S)+\Gamma_{\cd(V)}$ \`a un ordre au moins $(\cd(V)+u)T$. Ce qui pr\'ec\`ede  permet ``presque'' de d\'eduire que la section $D(\sigma)$ est nulle sur $\Gamma_{\cd(V)}+V$ pour tout op\'erateur diff\'erentiel $D$ d'ordre inf\'erieur \`a $\cd(V)T$. Autrement dit, au vu de la d\'efinition de l'ensemble $V_{\cd(V)-1}$ on a ``presque'' l'inclusion voulue avec multiplicit\'e au moins $T$. Il reste \`a expliquer le ``presque'' : \'etant donn\'ee une section $\mathfrak{s}\in H^0(X,\L^{\otimes D})$ nulle sur $\Gamma(S)+\Gamma_{\cd(V)}$ \`a un ordre au moins $(\cd(V)+u)T$, le raisonnement du paragraphe pr\'ec\'edent appliqu\'e \`a $\mathfrak{s}$ permet de conclure que $\Gamma_{\cd(V)}+V\subset Z(\mathfrak{s})$. Donnons nous maintenant une sous-vari\'et\'e (non-n\'ecessairement irr\'eductible) $Y$ de $X$ contenant $\Gamma(S)+\Gamma_{\cd(V)}$ et $\Gamma_{\cd(V)}+V$, et contenue dans $Z(\mathfrak{s})$ (on peut par exemple prendre $Y=Z(\mathfrak{s})$ muni de sa structure de sch\'ema r\'eduit). La vari\'et\'e $Y$ a \'et\'e choisie telle que pour tout op\'erateur diff\'erentiel $D$ d'ordre $1$, $D(\mathfrak{s})$ est une section bien d\'efinie de $H^0(Y,\L_{|Y}^{\otimes D})$. On peut donc it\'erer le raisonnement de la fa\c{c}on suivante : on r\'eapplique l'argument du paragraphe pr\'ec\'edent avec la section $s=D(\sigma)$. Notons que cette fois, la section est d\'efinie sur $Y$ et non plus sur $X$. De m\^eme la section $t_{g}^*(s)$ est d\'efinie sur $Y_g:=t_{-g}(Y)$ et le lieu des z\'eros $Z(t^*_g(s))$ doit \^etre vu comme \'etant le lieu des z\'eros dans $Y_g\subset X$. Ceci \'etant, la vari\'et\'e $Y_g$ contenant $V$, le raisonnement reste inchang\'e. On en d\'eduit que $D(\sigma)$ est nulle sur $\Gamma_{\cd(V)}+V$ et en appliquant ce proc\'ed\'e jusqu`\`a l'ordre $\cd(V)T$ on en d\'eduit le r\'esultat voulu.\hfill$\Box$

\medskip

\noindent Le lemme suivant donne une valeur de $u$ pour laquelle on peut appliquer la proposition \ref{l1} pr\'ec\'eden\-te.

\medskip

\begin{lemm}\label{rajout}Soient $S$ et $T$ deux entiers strictement positifs, $\L$ un fibr\'e en droites ample, $V$ une vari\'et\'e exceptionnelle de Seshadri de $\L$ relativement \`a $\Gamma(S)$ et $u$ la solution r\'eelle strictement comprise entre $0$ et $1$ de l'\'equation $x^{{\dim X}-1}(x+\codim V)=1$. On suppose de plus que $T>\alpha\left(S,\L^{\otimes D}\right)$ et qu'il existe une section non nulle $\sigma$ de $H^0(X,\L^{\otimes D})$ s'annulant \`a un ordre au moins $(\codim V+u)T$ sur $\Gamma(S)$. Sous ces hypoth\`eses on a 
\[uT>\e\left(\Gamma(S),\L^{\otimes D}\right).\]
\end{lemm}
\demo On note $d=\dim X$. C'est un simple calcul. On a par d\'efinition de $\alpha(S,\cdot)$~:
\[ \alpha\left(S,\L^{\otimes D}\right)^{d}=D^d\alpha\left(S,\L\right)^{d}=D^{d}\frac{\deg_{\L}(X)}{|\Gamma(S)|}=D^{d-1}\frac{\L^{d-1}\cdot Z(\sigma)}{|\Gamma(S)|}.\]
\noindent On en d\'eduit  
\begin{align*}
\e\left(\Gamma(S),\L^{\otimes D}\right)	&	\leq	\left(\frac{\L^{d-1}\cdot Z(\sigma)}{|\Gamma(S)|(\codim V+u)T}\right)^{\frac{1}{d-1}}D\text{ par le lemme \ref{l0}}\\
										&	<		\alpha\left(S,\L^{\otimes D}\right)\left(\frac{1}{u+\codim V}\right)^{\frac{1}{d-1}}\text{ car }T>\alpha\left(S,\L^{\otimes D}\right)\\
										&	<	Tu.
\end{align*}
\noindent \hfill$\Box$

\medskip

\begin{lemm}\label{l2}Avec les notations qui pr\'ec\`edent on a la suite d'inclusions
\[\forall k\in\left\{0,\ldots,\cd(V)-1\right\}\ \ \ \ \Gamma(a_{\cd(V)-k})+V_{k+1}\subset V_k\]
\noindent et chacune des inclusions avec multiplicit\'e au moins $T$. De plus le dernier ensemble de la tour d'inclusions (l'ensemble $V_0$) est strictement inclus dans $X$.
\end{lemm}
\demo Le fait que le dernier ensemble de la suite d'inclusions soit strictement contenu dans $X$ r\'esulte de ce que la section $\sigma$ n'est pas la section nulle. Pour les inclusions correspondant \`a $0\leq k\leq \cd(V)-2$, l'affirmation r\'esulte simplement de la d\'efinition des ensembles $V_k$. Enfin pour $k=\cd(V)-1$ il s'agit pr\'ecis\'ement de la proposition \ref{l1} pr\'ec\'edente.\hfill$\Box$

\medskip

\noindent En oubliant les multiplicit\'es on voit que l'on a construit, comme annonc\'e une tour d'ensembles emboit\'es, de dimensions comprises entre $\dim V$ et $\dim X-1$. De plus cette tour contient exactement $\codim V$ ensembles. Par le principe des tiroirs au moins une des inclusions de la tour correspond \`a une \'egalit\'e de dimension. Notons $r$ le plus petit entier v\'erifiant $\dim V_r=\dim V_{r+1}$. Le reste de la preuve est ensuite exactement calqu\'e sur la m\'ethode de Philippon, revisit\'ee par Nakamaye, et n'apporte rien de nouveau. Nous en donnons une preuve rapide, bas\'ee sur \cite{nak1}, au paragraphe suivant.

\section{Fin de la preuve du point 1. du th\'eor\`eme \ref{t2}}

\noindent Notons $W_r$ une composante irr\'eductible de $V_{r+1}$ de dimension maximale. C'est aussi une composante irr\'eductible isol\'ee de $V_r$ et mieux, $\Gamma(a_{\cd(V)-r})+W_r$ est une somme de composantes irr\'eductibles isol\'ees de $V_r$. En utilisant que l'inclusion est avec multiplicit\'e $T$ nous allons rappeler comment dans cette situation, on construit en suivant la preuve classique des lemmes de z\'eros, en coupant par des diviseurs num\'eriquement \'equivalents \`a $\L^{\otimes D}$ un cycle $C_{\codim W_r}$ dont le support contient  $\Gamma(a_{\cd(V)-r})+W_r$ avec multiplicit\'e $T$. Le th\'eor\`eme de B\'ezout nous donnera alors l'in\'egalit\'e attendue 
\[T^{\codim W_r}\deg_{\L^{\otimes D}}\left(\Gamma(a_{\cd(V)-r})+W_r\right)\leq\deg_{\L^{\otimes D}} X.\]

\medskip

Soit $C_1 = Z(\sigma)$.  Si $r = 0$ alors toute composante irr\'eductible de $\Gamma(a_{\cd(V)-r}) + W_r$ est une composante irr\'eductible de $C_1$ avec multiplicit\'e $T$ et il n'y a rien \`a faire. On suppose maintenant que $r \geq 1$ et on \'ecrit 
\[C_1 = C_1^\prime + \sum_{i=1}^k a_iZ_i,\]
\noindent o\`u les $Z_i$ sont les composantes irr\'eductibles de $C_1$ d'intersection non vide avec $\G$ et o\`u $C_1'$ est de support dans le diviseur \`a l'infini $X\backslash \G$.  

\medskip

\noindent Soit $Z_i$ une composante irr\'eductible de $C_1$ telle que $Z_i \cap \G$ est non-vide.

\begin{lemm} Il existe $x_i \in \Gamma_r$  et un op\'erateur diff\'erentiel $D_i$ d'ordre au plus $rT$ tel que la section $t_{x_i}^\ast(D_i(\sigma)) \in 
H^0(Z_i, t_{x_i}^\ast(\L))$ est bien d\'efinie et non nulle sur $Z_i$.  
\label{lu}
\end{lemm}
\demo Il s'agit du lemme 6 de \cite{nak1}.\hfill$\Box$

\medskip

\noindent Le cas des composantes irr\'eductibles $Z$ 
de $C_1$ contenues dans le diviseur \`a l'infini $X\backslash \G$ (\textit{i.e.} le cas des composantes irr\'eductibles de $C_1'$) est un peu plus d\'elicat. Le fibr\'e $\L$ \'etant ample, il existe un entier $m > 0$ tel que $\L^{\otimes mD}$ est tr\`es ample. Ainsi la classe $c_1(mD\L) \cap Z$ est repr\'esent\'ee par un cycle effectif $E_Z$ et on a 
\[c_1(\L^{\otimes D}) \cap Z = \frac{1}{m} E_Z.\]
\noindent On note $C_2'$ le cycle $\qq$-effectif ainsi obtenu : 
\[C_2'=\frac{1}{m}\sum_{Z\cap \G=\emptyset}E_Z.\]

\noindent Nous pouvons maintenant d\'efinir un cycle $C_2$ de codimension 2~: pour chacune des composantes $Z_i$ d'intersection non vide avec $\G$, le lemme \ref{lu} pr\'ec\'edent nous fournit un point $x_i\in \Gamma_r$ et un op\'erateur $D_i$ d'ordre au plus $rT$. On pose
\[C_2 = C_2^\prime + \sum_{i=1}^k a_i Z\left(t_{x_i}^\ast(D_i(\sigma))|Z_i\right).\]

\begin{lemm} On a $C_2 \equiv c_1(\L)^2$.
\label{lw}
\end{lemm}
\demo C'est le lemme 7 de \cite{nak1}.\hfill$\Box$

\medskip

\noindent Repartant de $C_2$ en lieu et place de $C_1$, nous pouvons it\'erer ce proc\'ed\'e jusqu'\`a atteindre un cycle $C_t$ de codimension $t=\codim W_r$ avec 
\[C_t \equiv c_1(\L)^t.\]
\noindent Par construction, toute composante irr\'eductible de dimension maximale de $V_r$ est une composante irr\'eductible de $C_t$ avec multiplicit\'e $T$. Ceci est donc en particulier vrai pour les composantes de $\Gamma(a_{\cd(V)-r})+W_r$ et le th\'eor\`eme de B\'ezout nous donne
\begin{equation}\label{varobs}
T^{\codim W_r}\deg_{\L}(\Gamma(a_{\cd(V)-r})+W_r)\leq\deg_{\L}(X)D^{\codim W}.
\end{equation}
\medskip
\noindent L'in\'egalit\'e attendue dans la conclusion du point 1. du th\'eor\`eme est donc obtenue, mais avec une vari\'et\'e $W_r$ et non pas un sous-groupe obstructeur comme annonc\'e. Pour cela introduisons $H=\stab(W_r)^0$ la composante connexe de l'identit\'e du stabilisateur $\stab(W_r\cap \G)$. Pour tout $x\in \Gamma(a_{\cd(V)-r})$, $x+\overline{H}$ est une composante irr\'eductible de $\bigcap_{y\in W_r\cap \G}t^*_y(V_r)$. Ainsi en repartant du cycle $C_t$ et en continuant la m\^eme proc\'edure en intersectant cette fois par des diviseurs de la forme $t^*_{y}(t^*_{x}(D(\sigma)))$  (au lieu simplement comme pr\'ec\'edemment des $t^*_{x}(D(\sigma))$) avec
$y \in W_r \cap \G$, $x \in \Gamma_r$ et $D$ op\'erateur diff\'erentiel d'ordre au plus $rT$, on aboutit \`a un cycle $C_h$ de codimension $h = \codim H$ contenant $x+\overline{H}$ avec multiplicit\'e $T$ pour tout $x\in\Gamma(a_{\cd(V)-r})$. Le th\'eor\`eme de B\'ezout donne~:
\begin{equation}\label{grobs}
T^{\codim H}\deg_{\L}(\Gamma(a_{\cd(V)-r})+\overline{H})\leq\deg_{\L}(X)D^{\codim H}.
\end{equation}
\noindent En remarquant que $\deg_{\L}(\Gamma(a_{i})+\overline{H})=\card\left(\Gamma(a_{i})+H/H\right)\deg_{\L}(\overline{H})$ ceci conclut.

\section{Preuve du point 2. du th\'eor\`eme \ref{t2}}

\noindent Pour prouver ce r\'esultat, il suffit de repartir de l'in\'egalit\'e (\ref{varobs}) que l'on a obtenue au paragraphe pr\'ec\'edent~: il existe un entier $r$ compris entre $0$ et $\codim V-1$ et il existe une vari\'et\'e $W_r$ de dimension inf\'erieure \`a $\dim X-1-r$ contenant un translat\'e de $V$ par un point de $\Gamma(a_1+\cdots+a_{\cd(V)-r-1})$ (quand $r=\cd(V)$ cette notation repr\'esente par convention l'ensemble $\{0\}$) telle que 
\begin{equation}\label{d1}
T^{\codim W_r}\deg_{\L^{\otimes D}}(\Gamma(a_{\cd(V)-r})+W_r)\leq\deg_{\L^{\otimes D}}(X)=\deg_{\L}(X)D^{\dim X},
\end{equation}
\noindent Supposons par l'absurde que $H$ soit de dimension $0$ autrement dit que $\stab(W_r\cap\G)$ est uniquement constitu\'e de points de torsion. Notons $b$ le nombre de composantes irr\'eductibles de $\Gamma(a_{\cd(V)-r})+W_r$. Une telle composante est de la forme $x+W_r$ avec $x\in\Gamma(a_{\cd(V)-r})$. Or
\begin{eqnarray*} 
\stab(W_r) & =			& \left\{x\in\G\ / x+W_r=W_r\right\}\\
			&\subset	& \left\{x\in\G\ / x+(W_r\cap\G)=(W_r\cap \G)\right\}=\stab(W_r\cap\G).
\end{eqnarray*}
Ainsi, si $x,y\in\Gamma(a_{\cd(V)-r})$, on constate que 
\[x+W_r=y+W_r\Rightarrow x-y\in \stab(W_r\cap\G)\cap\left(\Gamma(a_{\cd(V)-r})-\Gamma(a_{\cd(V)-r})\right).\]
\noindent Or ce dernier ensemble est inclus dans $\Gamma_{\mathbb{Z},\, \tors}.$ En particulier on a la majoration $b\geq \frac{|\Gamma(a_{\cd(V)-r})|}{|\Gamma_{\mathbb{Z},\, \tors}|}$. Revenant aux calculs de degr\'es, ceci entraine
\begin{eqnarray}
\deg_{\L^{\otimes D}}(\Gamma(a_{\cd(V)-r})+W_r)	&	\geq	&	 D^{\dim W_r}\frac{|\Gamma(a_{\cd(V)-r})|}{|\Gamma_{\mathbb{Z},\, \tors}|}\deg_{\L}(W_r) \\
												&	\geq	&	 D^{\dim W_r}\frac{|\Gamma(a_{\cd(V)-r})|}{|\Gamma_{\mathbb{Z},\, \tors}|}\label{d2}.
\end{eqnarray}												
\noindent En mettant ensemble les deux in\'egalit\'es (\ref{d1}) et (\ref{d2}), et en utilisant l'in\'egalit\'e v\'erifi\'ee par d\'efinition par $T$ (\textit{i.e.} $T\geq \alpha(S,\L^{\otimes D})$), on en d\'eduit que 
\begin{eqnarray*}
\left|\Gamma(a_{\cd(V)-r})\right| 	& \leq & |\Gamma_{\mathbb{Z},\, \tors}|\deg_{\L}(X)^{\frac{\dim W_r}{\dim X}}\left|\Gamma(S)\right|^{\frac{\codim W_r}{\dim X}}\\
									& \leq & |\Gamma_{\mathbb{Z},\, \tors}|(\deg_{\L} (X))\left|\Gamma(S)\right|^{\frac{\codim V}{\dim X}}.
\end{eqnarray*}

\noindent Le choix des $a_r$ permet de conclure par l'absurde.

\section{Preuve du th\'eor\`eme \ref{t4}}

\noindent Dans la suite nous noterons $j=\frac{1}{2}(1+\sqrt{5})$ le nombre d'or et $\equiv$ l'\'equivalence num\'erique. 

\medskip

\noindent \textbf{D\'emonstration du th\'eor\`eme \ref{t4}~: } 
\noindent Notons tout d'abord qu'il suffit de prouver le th\'eor\`eme pour $D=1$. En effet si le r\'esultat est acquis pour $D=1$, alors en prenant $D$ quelconque et en posant $\mathcal{M}=\L^{\otimes D}$, le r\'esultat appliqu\'e \`a $\mathcal{M}$ est exactement le r\'esultat voulu (on utilise ici la remarque 2 et le fait que $\deg_{\L^{\otimes D}}V=D^{\dim V}\deg_{\L}V$). On suppose donc d\'esormais que $D=1$. Si $C$ est une vari\'et\'e exceptionnelle de Seshadri du fibr\'e $\L$ relativement \`a $\Gamma(S)$ alors $C$ est n\'ecessairement de dimension $1$. En effet, sinon on aurait $C=X$ et le lemme \ref{l0} pr\'ec\'edent nous donnerait les in\'egalit\'es
\begin{align}
\left(\frac{\L^2}{|\Gamma(S)|}\right)^{\frac{1}{2}}& = \e(\Gamma(S),\L)\leq\frac{\L\cdot Z(\sigma)}{\sum_{x\in\Gamma(S)}\mult_x Z(\sigma)}\\
														& \leq \frac{\L^2}{|\Gamma(S)|jT}\text{ car $0\in \Gamma$ et $Z(\sigma)\equiv\L$}\label{or}\\ 
														& < \left(\frac{\L^2}{|\Gamma(S)|}\right)^{\frac{1}{2}}     \text{ par hypoth\`ese sur $T>\alpha(S,\L)=\left(\frac{\L^2}{|\Gamma(S)|}\right)^{\frac{1}{2}}$.}
\end{align}
Ceci est impossible. Par ailleurs, l'in\'egalit\'e $(\ref{or})$ nous indique que 
\[\e(\Gamma(S),\L) \leq \frac{\alpha(S,\L)^2}{jT} < j^{-1}T.\]
\noindent Utilisant que $jT=(1+j^{-1})T$, la proposition \ref{l1} appliqu\'ee avec $u=j^{-1}$ prouve que la section $\sigma$ s'annule le long de $x + C$ \`a un ordre sup\'erieur \`a $T$ pour tout $x \in \Gamma(a)$ (notons que pour pouvoir parler de $x+C$ nous utilisons ici l'hypoth\`ese faite sur la compactification). On en d\'eduit que $\L - T(\Gamma(a) + C)$ est un diviseur effectif et en intersectant avec $\L$, ceci produit l'in\'egalit\'e annonc\'ee.  

\medskip

\noindent Il reste \`a montrer que $C$ est l'adh\'erence d'un translat\'e de groupe alg\'ebrique. Notons pour cela que 
\[\stab(C)=\left\{x\in\G\ / \ x+C=C\right\}\subset\stab(C\cap\G)=\bigcap_{y\in C\cap\G}\left(y-(C\cap \G)\right).\]
\noindent Ainsi, si le groupe $\stab(C)$ est de dimension $1$, il en est de m\^eme du groupe $\stab(C\cap \G)$. En particulier $C\cap \G$ est un translat\'e d'un groupe alg\'ebrique de dimension $1$, dont l'adh\'erence est $C$. On peut donc supposer maintenant que $\stab(C)$ et m\^eme $\stab(C\cap \G)$ est de dimension $0$, autrement dit est uniquement constitu\'e de points de torsion. On a pour tout entier $k>0$ une inclusion de $H^0(A,kT(\Gamma(a)+C))$ dans $H^0(A,k\L)$. On en d\'eduit
\begin{equation}\label{eq1}
\L^2 \geq T^2(\Gamma(a)+C)^2 > \alpha(S,\L)^2(\Gamma(a)+C)^2.
\end{equation}
\noindent Notons $b$ le nombre de composantes irr\'eductibles de $\Gamma(a)+C$. Pour minorer $b$ nous allons majorer le nombre de $x_i,x_j \in \Gamma(S)$ tels que 
\[x_i + C = x_j + C.\]
\noindent Autrement dit, nous voulons majorer le nombre de $y \in \Gamma(S) - \Gamma(S)\subset \Gamma_{\mathbb{Z}}$ tels que
\[y + C = C.\]
Le stabilisateur de $C$ \'etant fini par hypoth\`ese, un tel $y$ est n\'ecessairement dans $\Gamma_{\mathbb{Z},\, \tors}$. On en d\'eduit une minoration du nombre de composantes irr\'eductibles~:  
\begin{equation}\label{eq100}
b \geq \frac{|\Gamma(a)|}{\left|\Gamma_{\mathbb{Z},\, \tors}\right|}.
\end{equation}
\noindent En injectant ceci dans l'in\'egalit\'e (\ref{eq1}) et en appliquant l'hypoth\`ese sur $a$, on obtient
\begin{equation}\label{eq101}
\L^2 > \L^2\frac{C^2 b^2}{|\Gamma(S)|} >\L^2C^2.
\end{equation}
\noindent Il suffit de noter que $C^2 \geq 1$ pour conclure par l'absurde. En effet, $\stab(C\cap \G)$ est de dimension $0$, donc $C\cap\G$ n'est pas un translat\'e de groupe alg\'ebrique et donc 
\[C^2 = (C-\eta_1) \cap (C-\eta_2) \supset {\Stab(C_{\G})}\]
o\`u $\eta_1$ et $\eta_2$ sont deux points g\'en\'eraux de $C\cap\G$.  

\section{Cas des d\'erivations le long d'un sous-espace de dimension $1$ dans une surface\label{observation}}

\noindent On peut se demander ce qui se passe dans le th\'eor\`eme \ref{t4} quand le sous-espace $\EV \subset T_0(\G)$ est de dimension 1. C'est l'objet du th\'eor\`eme \ref{t5} ci-dessous. 

\medskip

\noindent Notons tout d'abord que la complexit\'e apport\'ee \`a
la preuve ne se trouve pas du cot\'e des d\'eriva\-tions~: si $D_v$ est la d\'erivation correspondant \`a un vecteur $0 \neq v \in  \EV$ on a pour toute courbe $C$, 
\[C \subset Z(\sigma)\Rightarrow \mult_C(D_v(\sigma)) < \mult_C(\sigma),\] 
sauf si la courbe $C_{\G}$ est un translat\'e d'un sous groupe dont l'espace tangent \`a l'origine est pr\'ecis\'ement $\EV$.

\medskip

\noindent Le probl\`eme vient de ce que une fois que l'on consid\`ere la multiplicit\'e le long d'un sous-espace propre, la notion de sous-vari\'et\'e exceptionnelle de Seshadri n'a plus de sens. En effet, pour d\'efinir la constante de Seshadri on utilise le fait que l'on peut se donner une r\'esolution simultan\'ee de 
tous les id\'eaux de la forme $m_x^t$ o\`u $t$ est un entier, $x \in X$, et $m_x \subset \oo_X$ est l'id\'eal maximal associ\'e au point $x$. Cette r\'esolution permet de faire varier le diviseur exceptionnel d'une mani\`ere continue et la d\'efinition de la constante de Seshadri exploite ce fait~: en particulier, si $\pi: \tilde{X} \ra X$ est l'\'eclatement
de $X$ en $x$, de diviseur exceptionnel $E$,
 alors la constante de Seshadri d'un fibr\'e en droites ample $\L$
sur $X$ peut \^etre d\'efinie par
\[\e(x,\L) = \sup_{\alpha \in {\bf R}}\{\pi^\ast(\L)(-\alpha E) \,\,
\mbox{est nef}\}.\]

\medskip

\noindent Quand il s'agit de l'ordre d'annulation d'une section le long d'un sous-espace strict $\EV$, nous pouvons tout de m\^eme dire quelque chose :
 pour tout entier $k$, soit 
$\pi_k: X_k \ra X$ une r\'esolution lisse
du faisceau d'id\'eaux
\[\mathcal{I}_k = \bigcap_{x \in \Gamma(S)} \mathcal{I}_{x,k}\]
o\`u $\mathcal{I}_{x,k}$ est engendr\'e par les fonctions r\'eguli\`eres $f$, 
d\'efinies dans un voisinage de $x$, telles que $f$ s'annule le long 
de $\EV$ \`a un ordre au moins $k$. Autrement dit, on fait \'eclater
simultan\'ement chaque id\'eal $\mathcal{I}_{x,k}$ et puis on prend une
r\'esolution de cette vari\'et\'e: en particulier, l'application
$\pi_k: X_k \ra X$ est birationelle et un isomorphisme hors des points
$x \in \Gamma(S)$.
On note $E_k$ le diviseur exceptionnel de $\pi_k : X_k \ra X$.

\medskip

\defi Nous pouvons d\'efinir un \textit{analogue} de la constante de Seshadri par
\[\e_{\EV}(\Gamma(S),\L) = \sup_{k \in \N} \left\{\pi_k^\ast(\L)(-E_k) \ \text{ est nef}\right\}.\]
\noindent On dira que $\e_{\EV}(\Gamma(S),\L)$ est la \textit{constante de Seshadri de $\L$ relativement \`a $\Gamma(S)$ et $\EV$}. Par ailleurs on appelle \textit{sous-vari\'et\'e exceptionnelle de Seshadri relativement \`a $\Gamma(S)$ et $\EV$} toute vari\'et\'e  $V\subset X$ telle que 
\[(\pi^\ast(\L)(-E_{\e_\EV(\Gamma(S),\L)+1}))^{\dim(V)} 
\cdot \tilde{V} < 0\]
o\`u $\tilde{V}$ est la transform\'ee stricte de $V$ dans $X_k$. 

\medskip 

\remarque Il s'agit bien d'un analogue et pas d'une g\'en\'eralisation : quand $\EV=T_0(\G)$, la constante $\e_{\EV}(\Gamma(S),\L)$ est toujours un
entier tandis que $\e(\Gamma(S),\L)$ peut \^etre rationnelle ou m\^eme \textit{a priori} irrationnelle. Par contre on voit sur les d\'efinitions que l'on a toujours l'in\'egalit\'e $\e_{\EV}(\Gamma(S),\L)\leq \e(\Gamma(S),\L)$,
avec \'egalit\'e si et seulement si $\e(\Gamma(S),\L)$ est un entier.

\medskip

\remarque Cette notion de constante de Seshadri relative \`a $\EV$ n'est pas, en g\'en\'eral, homog\`ene.  Autrement dit il n'est pas vrai en g\'en\'eral que
\[\e_{\EV}(\Gamma(S),\L^{\otimes n}) =n\e_{\EV}(\Gamma(S),\L).\]
Cette constante $\e_{\EV}(\Gamma(S),\L^{\otimes n})$ peut grossir comme un polyn\^ome en $n$.

\medskip

\remarque Il faut en g\'en\'eral remplacer $\L$ par $\L^{\otimes {d}}$ avec $d\gg 0$ pour que
$\e_{\EV}(\Gamma(S),\L)$ soit non nulle.  Ce probl\`eme n'existe \'evidemment
pas avec la constante de Seshadri $\e(\Gamma(S),\L)$ qui peut
\^etre petite et positive (alors que dans le cas de notre constante relative, une fois que
$\e_{\EV}(\Gamma(S),\L) < 1$ on a n\'ecessairement
$\e_{\EV}(\Gamma(S),\L) = 0$).  

\medskip

\noindent Nous n'utiliserons pas les remarques pr\'ec\'edentes dans la suite.

\medskip

\begin{lemm}\label{lf}  On a l'in\'egalit\'e
\[\e_\EV(\Gamma(S),\L) \leq \alpha(S,\EV,\L).\]
\end{lemm}
\demo Supposons le r\'esultat faux et raisonnons par l'absurde. Dans ce cas il existe un entier $\beta > \alpha(S,\EV,\L)$ tel que
le fibr\'e $\pi_\beta^\ast \L(-E_\beta)$ soit nef.  Par cons\'equent on a
\[c_1(\pi_\beta^\ast \L(-E_\beta))^{2} \geq 0.\]

\noindent Soit alors $x$ un point de $X$ et $\pi: Y \ra X$ une r\'esolution de l'id\'eal $\mathcal{I}_{x,\beta}$ de diviseur exceptionnel $E_x$. On a
\begin{eqnarray*}
c_1(\pi_\beta^\ast \L(-E_\beta))^{2} & =&   c_1(\pi_\beta^\ast({\cal L}))^2 +
c_1({\mathcal O}_Y(-E_\beta))^2  \\
& = &
c_1(\L)^{2} + |\Gamma(S)|
c_1(\mathcal{O}_Y(-E_x))^{2}.
\end{eqnarray*}
D'autre part 
\[c_1(\mathcal{O}_Y(-E_x))^{2} = -\beta.\]


\medskip

\noindent On en d\'eduit donc 
\[\deg_{\L}(X)=c_1(\L)^{2} \geq |\Gamma(S)|\beta.\]
\noindent Le sous-espace $\EV$ \'etant de dimension $1$, cette derni\`ere in\'egalit\'e contredit la d\'efi\-ni\-tion de $\alpha(S,\EV,\L)$.\hfill$\Box$

\medskip

\noindent Nous pouvons maintenant \'enoncer le r\'esultat principal de ce paragraphe~: 

\begin{theo}\label{t5} On suppose que $\dim X=2$ et 
que $\Gamma_{\mathbb{Z}}$ n'est pas de torsion. Soient $\EV \subset T_0(\G)$
un sous-espace de dimension 1, et $S$, $a$ deux entiers strictement positifs tels 
que $|\Gamma(S)|\geq |\Gamma(a)|\geq|\Gamma_{\mathbb{Z},\textnormal{tors}}|
\cdot|\Gamma(S)|^{\frac{1}{2}}$. Soient $T>\alpha(S,\EV,\L^{\otimes D})$ un entier 
et $\sigma \in H^0(X,\L^{\otimes D})$ une section non nulle s'annulant le long de $\EV$ \`a un ordre 
au moins $T + \e_\EV(\Gamma(S),\L)$ sur $\Gamma(S+a)$. 
Toute sous-vari\'et\'e exceptionnelle de Seshadri de 
$\L$ relativement \`a $\Gamma(S)$ et $\EV$ est l'adh\'erence d'un translat\'e 
de sous-groupe alg\'ebrique de dimension $1$.  Soit $E$ l'une de
ces courbes.
\begin{enumerate}
\item Si $T_0(E) \neq \EV$ alors
$\sigma$ s'annule sur
$\Gamma(a) + E$ \`a un ordre le long de $\EV$ au moins $T$ et on a
\[ T\cdot\card(\Gamma(a) + E/E)\deg_{\L}(E) \leq D\deg_{\L}(X).\]
\item Si $T_0(E) = \EV$ alors $\sigma$ s'annule sur $\Gamma(a) + E$ et on a
\[ \card(\Gamma(a) + E/E)\deg_{\L}(E) \leq D\deg_{\L}(X).\]
\end{enumerate}
\end{theo}

\medskip

\remarque Comme pour le th\'eor\`eme \ref{t4}, le lemme \ref{lf} pr\'ec\'edent nous assure que la condition d'annulation demand\'ee  : $T + \e_\EV(\Gamma(S),\L)\leq T+\alpha(S,\EV,\L)<2T$ est plus faible que la condition classique $2T+1.$  Par ailleurs  comme l'indique la preuve, l'hypoth\`ese sur la taille de $a$ par rapport \`a $S$ (qui est la m\^eme hypoth\`ese que celle faite dans le th\'eor\`eme \ref{t4}) n'intervient que pour prouver que la courbe exceptionnelle est l'adh\'erence d'un groupe al\'egbrique, autrement dit pour prouver que le groupe obstructeur qui intervient est \textit{non nul}. Cette hypoth\`ese n'intervient pas dans la preuve de l'in\'egalit\'e num\'erique sur le degr\'e de la courbe obstructrice.

\medskip

\noindent La preuve est calqu\'ee sur celle du th\'eor\`eme \ref{t4}. On constate tout d'abord que l'on peut l\`a encore supposer $D=1$. Ce que l'on fait d\'esormais. Ensuite  il faut adapter la proposition \ref{l1} \`a notre nouveau cadre (d\'erivation le long d'un sous-espace de dimension $1$)~: c'est l'objet du lemme qui suit. Soit $C \subset X$ une courbe qui est une sous-vari\'et\'e
exceptionnelle de Seshadri relativement \`a $\Gamma(S)$ et $\EV$. Il y a deux possibilit\'es~: soit $C$ est un 
translat\'e d'un sous--groupe $H \subset G$ avec $T_0(H) = \EV$ ; soit ce n'est
pas le cas.  La raison pour laquelle le cas des translat\'es est particulier est
que si la section $\sigma \in H^0(X,\L)$ s'annule sur un tel translat\'e
alors son ordre d'annulation le long de $\EV$ est infini.  

\medskip

\begin{lemm} Soient $T,S,a$ trois entiers strictement positifs, $\L$ un fibr\'e en droites ample et $\sigma$ une section non-nulle de $H^0(X,\L)$. Supposons que $\sigma$ s'annule le long de $\EV$ \`a un ordre $\geq T + \e_\EV(\Gamma(S),\L)$ sur $\Gamma(S+a)$ et notons $C$ une courbe exceptionnelle de Seshadri de $\L$ relativement \`a $\Gamma(S)$ et $\EV$. Alors $\sigma$ s'annule le long
de $\EV$ \`a un ordre au moins $T$ sur $\Gamma(a) + C$.
\label{ll}
\end{lemm}
\demo Montrons d'abord que $\sigma$ s'annule sur $\Gamma(a) + C$. Posons $\e = \e_\EV(\Gamma(S),\L)$. Soit $g\in\Gamma(a)$. Par hypoth\`ese on sait que $s=t^\ast_g(\sigma)$ est nulle sur $\Gamma(S)$ le long de $\EV$ \`a un ordre au moins $\e+1$. Autrement dit, le diviseur $\pi^\ast(Z(s)) - E_{\e+1}$ est effectif et par cons\'equent en intersectant avec $\tilde{C}$, on en d\'eduit que cette intersection est positive, sauf si $\tilde{C}$ est une composante de ce diviseur. Or par la d\'efinition de $\e_\EV(\Gamma(S),\L)$ et le fait que
$C$ est une courbe exceptionnelle de Seshadri de $\L$ relativement \`a
$\Gamma(S)$ et $\EV$ on a,
\[\left(\pi^\ast(Z(s)) - E_{\epsilon+1}\right)\cdot\tilde{C}\equiv\pi_{\e+1}^\ast \L(- E_{\e +1}) \cdot \tilde{C} < 0.\]
\noindent  Ceci prouve donc que $\tilde{C}$ est une composante de $\pi^\ast(Z(s)$ et donc que $s$ s'annule
sur $C$.  Autrement dit $\sigma$ s'annule sur $\Gamma(a)+C$.

\medskip

\noindent Soient $v\in\EV$ un vecteur non nul de d\'erivation correspondante $D_{\EV}$, $i<T$ un entier et $x \in \Gamma(a)$. Supposons que
$D^i_\EV(\sigma)|(x+C)$ est une section bien d\'efinie de $H^0(x+C,\L)$.  
Localement la section $D^i_\EV(\sigma)$ s'identifie avec une fonction r\'eguli\`ere $f$ qui
s'annule le long de $\EV$
 sur $x+\Gamma(S)$ \`a un ordre au moins $\e_\EV(\Gamma(S),\L)+1$.
La courbe $C$ \'etant une courbe exceptionnelle de Seshadri de $\L$ relativement $\EV$ et $\Gamma(S)$, on en d\'eduit que $f$ s'annule identiquement sur $x+C$.\hfill$\Box$

\medskip

\noindent Nous pouvons maintenant donner la preuve du th\'eor\`eme \ref{t5}~: 
\demo
Notons $\e=\e_\EV(\Gamma(S),\L)$. L'existence d'une section $\sigma$ nulle sur $\Gamma(S)$ \`a un ordre au moins $T>\alpha(S,\EV,\L)$ entra\^ine que $\e_\EV(\Gamma(S),\L) < \alpha$. Or ceci implique que $\pi^\ast(\L(-E_{\e}))^2 > 0$ et donc que $X$ n'est pas une sous-vari\'et\'e exceptionnelle de Seshadri de $\L$ relativement \`a $\Gamma(S)$ et $\EV$.  Ainsi toute vari\'et\'e exceptionnelle est une courbe. Soit $C \subset X$ une courbe
qui est une sous-vari\'et\'e exceptionnelle de Seshadri pour $\L$ relativement \`a $\Gamma(S)$ et $\EV$.  Si $C \cap \G$ est un sous-groupe de $\G$
avec $T_0(C) = \EV$ il n'y a rien \`a montrer~:
d'apr\`es le lemme \ref{ll}
la section $\sigma$ s'annule sur $\Gamma(a) + C$
d'o\`u l'in\'egalit\'e.  Supposons donc que $C \cap \G$ n'est pas un translat\'e
d'un sous-groupe dont l'espace tangent \`a l'origine est $\EV$.  
Dans ce cas, si $\eta \in C$ est un point g\'en\'eral alors $T_\eta(C) \neq
\EV$ et par cons\'equent la d\'erivation $D_\EV$ est transverse \`a $C$ et \`a toutes
ses translat\'ees.  Par le lemme \ref{ll}, on voit que l'hypoth\`ese sur l'ordre d'annulation de $\sigma$ le long
de $\EV$ sur $\Gamma(a+S)$ garantit que, pour tout entier $i<T$, la section $D_\EV^i(\sigma)$ s'annule sur
$\Gamma(a) + C$.  Ainsi le diviseur
\[Z(\sigma) - T\sum_{x \in \Gamma(a)} (x+C)\]
est effectif.  On en d\'eduit l'in\'egalit\'e annonc\'ee dans le th\'eor\`eme \ref{t5}.
Quant au fait que dans ce cas $C \cap \G$ est un sous-groupe, la preuve est la m\^eme que celle (du m\^eme fait) du th\'eor\`eme \ref{t4}. C'est dans cette partie uniquement qu'est utilis\'ee l'hypoth\`ese faite sur la taille de $a$ relativement \`a celle de $S$.\hfill$\Box$



\begin{thebibliography}{1}

\bibitem{bauer}
T.~Bauer.
\newblock Seshadri constants on algebraic surfaces.
\newblock {\em Math. Ann.}, 313(3):547--583, 1999.

\bibitem{cp}
F.~Campana and T.~Peternell.
\newblock Algebraicity of the ample cone of projective varieties.
\newblock {\em J. Reine Angew. Math.}, 407:160--166, 1990.

\bibitem{fulton}
W.~Fulton.
\newblock {\em Intersection Theory}.
\newblock Springer, second edition, 1998.

\bibitem{lazarsfeld}
R.~Lazarsfeld.
\newblock {\em Positivity in algebraic geometry. {I}}, volume~48 of {\em
  Ergebnisse der Mathematik und ihrer Grenzgebiete. 3. Folge. A Series of
  Modern Surveys in Mathematics}.
\newblock Springer-Verlag, Berlin, 2004.

\bibitem{nak1}
M.~Nakamaye.
\newblock Multiplicity estimates on commutative algebraic groups.
\newblock \`A paraitre \`a Crelle's Journal, 2006.

\bibitem{philippon}
P.~Philippon.
\newblock Lemmes de z\'eros dans les groupes alg\'ebriques commutatifs.
\newblock {\em Bull. Soc. Math. France}, 114(3):355--383, 1986.

\bibitem{szpiro}
L.~Szpiro.
\newblock Sur les solutions d'un syst\`eme d'\'equations polynomiales sur une
  vari\'et\'e ab\'elienne (d'apr\`es {G}. {F}altings et {P}. {V}ojta).
\newblock {\em Ast\'erisque}, (189-190):Exp.\ No.\ 729, 429--446, 1990.
\newblock S\'eminaire Bourbaki, Vol.\ 1989/90.

\bibitem{serre}
M.~Waldschmidt.
\newblock Nombres transcendants et groupes alg\'ebri\-ques.
\newblock {\em Ast\'erisque}, (69-70):218, 1987.
\newblock Avec des appendices de Daniel Bertrand et Jean-Pierre Serre.

\end{thebibliography}
\end{document}